\documentclass[12pt]{amsart}
\usepackage{amssymb, amsthm, enumerate}
\usepackage[utf8]{inputenc}
\usepackage{graphicx}
\usepackage{indentfirst}
\newtheorem{definition}{Definition}
\newtheorem{theorem}{Theorem}

\newtheorem{lemma}{Lemma}

\newtheorem{corollary}{Corollary}
\newenvironment{Pf}{ Proof.}{\(\square\)}

\begin{document}

\title[On the generalization of Erd\H{o}s-Vincze's theorem...]{On the generalization of Erd\H{o}s-Vincze's theorem about the approximation of regular triangles by polyellipses in the plane}
\author{Cs. Vincze, Z. Kovács and Zs. F. Csorvássy}
\footnotetext[2]{{\bf MR subject classification:} 51M04}
\footnotetext[3]{Cs. Vincze is supported by the EFOP-3.6.1-16-2016-00022 project. The project is
co-financed by the European Union and the European Social Fund.}
\address{Institute of Mathematics, University of Debrecen, P.O.Box 12, H-4010 Debrecen, Hungary}
\email{csvincze@science.unideb.hu}
\address{Institute of Mathematics and Informatics, University of Nyíregyháza, H-4400 Nyíregyháza, \phantom{bl}Sóstói út 31./B, Hungary}
\email{kovacs.zoltan@nye.hu}
\address{Debrecen Reformed Theological University, H-4026 Debrecen, Kálvin tér 16, Hungary}
\email{csorvassy.zs.f@gmail.com}
\begin{abstract}
A polyellipse is a curve in the Euclidean plane all of whose points have the same sum of distances from finitely many given points (focuses). The classical version of Erdős-Vincze's theorem states that regular triangles can not be presented as the Hausdorff limit of polyellipses even if the number of the focuses can be arbitrary large. In other words the topological closure of the set of polyellipses with respect to the  Hausdorff distance does not contain any regular triangle and we have a  negative answer to the problem posed by E. Vázsonyi (Weissfeld) about the approximation of closed convex plane curves by polyellipses. It is the additive version of the approximation of simple closed plane curves by polynomial lemniscates all of whose points have the same product of distances from finitely many given points (focuses). Here we are going to generalize the classical version of Erdős-Vincze's theorem for regular polygons in the plane. We will conclude that the error of the approximation tends to zero as the number of the vertices of the regular polygon tends to the infinity. The decreasing tendency of the approximation error gives the idea to construct curves in the topological closure of the set of polyellipses. If we use integration to compute the average distance of a point from a given (focal) set in the plane then the curves all of whose points have the same average distance from the focal set can be given as the Hausdorff limit of polyellipses corresponding to partial sums. 
\end{abstract}
\maketitle

\section{Introduction}

Polyellipses in the plane belong to the more general idea of generalized conics \cite{VA_1}, \cite{VA_2}, \cite{VA_3}, see also \cite{GS}. They are subsets in the plane all of whose points have the same average distance from a given set of points (focal set). The level set of the function measuring the arithmetic mean of Euclidean distances from the elements of a finite pointset is one of the most important discrete cases. Curves given by equation 
\begin{equation}
\label{poly}
\sum_{i=1}^m d(X, F_i)=c \ \ \Leftrightarrow\ \ \frac{\sum_{i=1}^m d(X, F_i)}{m}=c_m \ \ (c_m:=c/m)
\end{equation}
are called polyellipses with focuses $F_1, \ldots, F_m$, where $d$ means the Euclidean distance in the plane. They are the additive version of lemniscates all of whose points have the same geometric mean of distances (i.e. their product is constant). Polyellipses appear in optimization problems in a natural way  \cite{MF}. The characterization of the minimizer of a function measuring the sum of distances from finitely many given points is due to E. Vázsonyi (Weissfeld) \cite{W}. He also posed the problem of the approximation of closed convex plane curves with polyellipses. P. Erdős and I. Vincze \cite{ErdosVincze} proved that it is impossible in general because regular triangles can not be presented as the Hausdorff limit of polyellipses even if the number of the focuses can be arbitrary large. The proof of the classical version of Erdős-Vincze's theorem can be also found in \cite{VinczeVarga}. The aim of the present paper is to generalize the theorem for regular polygons. Although a more general theorem can be found in P. Erdős and I. Vincze \cite{ErdosVincze1} stating that the limit shape of a sequence of polyellipses may have only one single straight segment, the high symmetry of regular polygons allows us to follow a special argument based on the estimation of the curvature of polyellipses with high symmetry. On the other hand we conclude that the error of the approximation is tending to zero as the number of the vertices of the regular polygon tends to the infinity. The decreasing tendency of the approximation error gives the idea to construct curves in the topological closure of the set of polyellipses. If we use integration to compute the average distance of a point from a given (focal) set in the plane then the curves all of whose points have the same average distance from the focal set can be given as the Hausdorff limit of polyellipses corresponding to partial sums. The idea can be found in \cite{ErdosVincze1} but it was formulated by some other authors as well, see e.g. \cite{GS}. 

Let $R>0$ be a positive real number. The \emph{parallel body} of a set $K$ with radius $R$ is the union of the closed disks with radius $R$ centered at the points of $K$. The infimum of the positive numbers such that $L$ is a subset of the parallel body of $K$ with radius $R$ and vice versa is called the \emph{Hausdorff distance} of $K$ and $L$. It is well-known that the Hausdorff metric makes the family of non\-empty closed and bounded (i.e. compact) subsets in the plane a complete metric space. Another possible characterization of the Hausdorff distance between compact subsets in the plane can be given in terms of distances between the points of the sets: if
$$\max_{X \in K} d(X, L):= \max_{X \in K}  \min_{Y \in L}d(X,Y) \ \ \textrm{and}\ \ \max_{Y \in L} d(K,Y):= \max_{Y \in L}  \min_{X \in K}d(X,Y)$$
then the Hausdorff distance of $K$ and $L$ is
\begin{equation}
\label{minimax}
h(A,B)= \max \{ \max_{X \in K} d(X, L), \max_{Y \in L} d(K,Y) \}.
\end{equation}

\section{Polyellipses in the Euclidean plane}

\begin{definition}
Let $F_1,\ldots,F_m$ be not necessarily different points in the plane and consider the function 
$$F(X):=\sum_{i=1}^{m}d(X,F_i)$$
measuring the sum of distances of $X$ from $F_1,\ldots,F_m$. The set given by equation $F(X)=c$ is called a polyellipse with focuses $F_1,\ldots,F_m$. The multiplicity of the focal point $F_i$ $(i=1, \ldots, m)$ is the number of its appearances in the sum.  
\end{definition}

\subsection{A Maple presentation I}

\begin{verbatim}
with(plottools):
with(plots):
#The coordinates of the focal points#
x:=[0.,1.,1.,0.,.5]: 
y:=[0.,0.,1.,1.,1.5]:
m:=nops(x):
F:=proc(u,v)
local i,partial:
partial:=0.:
for i to m do 
partial:=partial+sqrt((u-x[i])^2+(v-y[i])^2) 
end do 
end proc:
#The list of the level rates#
c:=[3.75,4,4.5,5,5.7]:
graf1:=pointplot(zip((s,t)->[s,t],x,y),symbol=solidcircle, symbolsize=20,
color=black):
graf2:=
seq(implicitplot(F(u,v)=i,u=-2..2,v=-2..2,numpoints=10000,color=black),i=c):
display(graf1,graf2,scaling=constrained,axes=none,linestyle=dot);
\end{verbatim}

\begin{figure}
\centering
\includegraphics[scale=0.3]{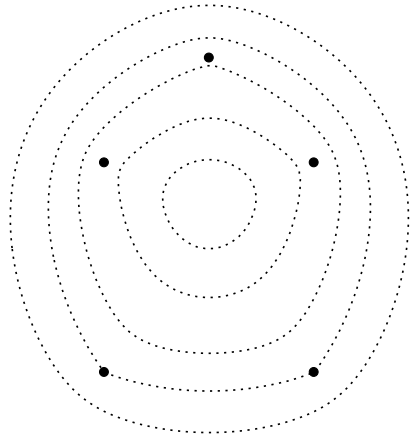}
\caption{A Maple presentation I.}
\end{figure}

\subsection{Weissfeld's theorems}
Using the triangle inequality it is easy to prove that $F$ is a convex function; moreover if the focuses are not collinear it is strictly convex. It is also clear that $F$ is differentiable at each point different from the focuses and
\begin{equation}
\label{firstder}
D_vF(X)=\langle v, N \rangle,\ \ \textrm{where}\ \ N:=-\sum_{i=1}^m \frac{1}{d(X,F_i)}\stackrel{\longrightarrow}{XF_i}
\end{equation}
is the opposite vector of the sum of unit vectors from $X$ to the focal points, respectively. Using that the vanishing of the first order derivatives is a sufficient and necessary condition for a point to be the minimizer of a convex function we have the first characterization theorem due to E. Vázsonyi (Weissfeld).

\begin{theorem} {\emph{\cite{W}}} Suppose that $X$ is different from the focal points; it is a minimizer of the function $F$ 
if and only if
$$N={\bf 0}\ \ \Leftrightarrow\ \ \sum_{i=1}^m \frac{1}{d(X,F_i)}\stackrel{\longrightarrow}{XF_i}={\bf 0}.$$
\end{theorem}

\noindent
A more subtle but standard convex analysis shows that if $X=F_i$ is one of the focal points  then
$$D_v^+F(F_i)=k_i\|v\|+\langle v, N_i\rangle,$$
where $D_v^+F(F_i)$ denotes the one-sided directional derivative of $F$ at $F_i$ into direction $v$ (the one-sided directional derivatives always exist in case of a convex function),
$$N_i:=- \sum_{F_j\neq F_i} \frac{1}{d(F_i,F_j)} \stackrel{\longrightarrow}{F_iF_j}$$
and $k_i$ is the multiplicity of the focal point $F_i$. Using that $D_v^+F(F_i)\geq 0$ for any direction $v$ 
is a sufficient and necessary condition for $F_i$ to be the minimizer of the function $F$, the Cauchy-Buniakovsky-Schwarz inequality gives the second characterization theorem due to E. Vázsonyi (Weissfeld).

\begin{theorem} {\emph{\cite{W}}} The focal point $X=F_i$ $(i=1, \ldots, m)$ is a minimizer of the function $F$ if and only if
$$\|N_i\|\leq k_i \ \ \Leftrightarrow\ \  \| \sum_{F_j\neq F_i} \frac{1}{d(F_i, F_j)} \stackrel{\longrightarrow}{F_iF_j} \|\leq k_i.$$
\end{theorem}

\noindent
In terms of coordinates we have the following formulas:

\begin{equation}
\label{gradient}
D_1F(x,y)=\sum_{i=1}^m\frac{x-x_i}{\sqrt{(x-x_i)^2+(y-y_i)^2}},
\end{equation}
$$D_2F(x,y)=\sum_{i=1}^m\frac{y-y_i}{\sqrt{(x-x_i)^2+(y-y_i)^2}}$$
provided that $X(x,y)$ is different from the focuses $F_i(x_i,y_i)$, where $i=1, \ldots, m$. The second order partial derivatives are
$$D_1D_1F (x,y)=\ \sum_{i=1}^m\frac{(y-y_i)^2}{\left((x-x_i)^2+(y-y_i)^2\right)^{3/2}},$$
\begin{equation}
\label{secondder}
D_2D_2F (x,y)=\ \sum_{i=1}^m\frac{(x-x_i)^2}{\left((x-x_i)^2+(y-y_i)^2\right)^{3/2}},
\end{equation}
$$D_1D_2F (x,y)=-\sum_{i=1}^m\frac{(x-x_i)(y-y_i)}{\left((x-x_i)^2+(y-y_i)^2\right)^{3/2}}.$$

\subsection{The main theorem} The classical version of Erdős-Vincze's theorem states that regular triangles can not be presented as the Hausdorff limit of polyellipses even if the number of the focuses can be arbitrary large and we have a  negative answer to the problem posed by E. Vázsonyi about the approximation of closed convex planar curves by polyellipses. It is the additive version of the approximation of simple closed planar curves by polynomial lemniscates all of whose points have the same product of distances from finitely many given points\footnote{The approximating process uses the partial sums of the Taylor expansion of holomorphic functions and the focuses correspond to the complex roots of the polynomials.}  (focuses). Here we are going to generalize the classical version of Erdős-Vincze's theorem for any regular polygon in the plane as follows.

\begin{theorem}
A regular $p$ - gon $(p\geq 3)$ in the plane can not be presented as the Hausdorff limit of polyellipses even if the number of the focuses can be arbitrary large. 
\end{theorem}

\section{The proof of Theorem 3}

Let $P$ be a regular $p$ - gon with vertices $P_1, \ldots, P_p$ inscribed in the unit circle centered at the point $O$ and suppose, in contrary, that there is a sequence $E_n$ of polyellipses tending to $P$. 

\subsection{The first step: the reformulation of the problem by circumscribed polyellipses.} Let $\varepsilon >0$ be an arbitrarily small positive real number and suppose that $n$ is large enough for $E_n$ to be contained in the ring of the circles $C_{-\varepsilon}$ and $C_{\varepsilon}$ around $O$ with radiuses $1-\varepsilon$ and $1+\varepsilon$. If $P_{-\varepsilon}$ is the regular $p$ - gon inscribed in $C_{-\varepsilon}$ then $E_n$ is a polyellipse around $P_{-\varepsilon}$. On the other hand 
$$h(P_{-\varepsilon}, E_n)\leq h(P_{-\varepsilon}, P)+h(P, E_n)\leq \varepsilon+h(C_{-\varepsilon},C_{\varepsilon})= 3\varepsilon.$$
Using a central similarity with center $O$ and coefficient $\frac{1}{1-\varepsilon}$ we have that $P_{-\varepsilon} \mapsto P_{-\varepsilon}'=P$ and $E_n \mapsto E_n'$, where $E_n'$ is a polyellipse around $P$ such that
$$h(P, E_n')=h(P_{-\varepsilon}',E_n')=\frac{1}{1-\varepsilon}h(P_{-\varepsilon}, E_n)\leq \frac{3\varepsilon}{1-\varepsilon}.$$
Taking $\varepsilon \to 0$ we have a sequence of circumscribed polyellipses tending to $P$. 

\subsubsection{Summary}
From now on we suppose that $E_n$ is a sequence of circumscribed polyellipses tending to $P$. 

\subsection{The second step: a symmetrization process.} Following the basic idea of \cite{ErdosVincze} we apply a symmetrization process to the circumscribed polyellipses of $P$ without increasing the Hausdorff distances. Let $F_1, \ldots, F_m$ be the focuses of a polyellipse $E$ around $P$ defined by the formula
$$\sum_{i=1}^{m}d(X,F_i)=c.$$
Consider the polyellipse $E_{sym}$ passing through the vertices $P_1, \ldots, P_p$ of $P$ such that its focal set is
\begin{equation}
\label{symmetry}
G:=\{f(F_i)|i=1,\ldots, m\ \textrm{and}\  f \in H\},
\end{equation}
where $H$ denotes the symmetry group of $P$. The equation of $E_{sym}$ is
$$\sum_{f \in H}\sum_{i=1}^{m}d(X,f(F_i))=c',$$
where
$$ c':=\sum_{f \in H}\sum_{i=1}^{m}d(P_1,f(F_i))=\sum_{f \in H}\sum_{i=1}^{m}d(P_2,f(F_i))= \ldots =\sum_{f \in H}\sum_{i=1}^{m}d(P_p,f(F_i))$$
because of the invariance of the vertices and the focal set under the group $H$. Note that $f(E)$ is a polyellipse around $P$ with focuses $f(F_1)$, $\ldots$, $f(F_m)$. It is defined by the equation
 \begin{equation}
\label{union}
\sum_{i=1}^{m}d(X,f(F_i))=c\ \ \textrm{and}\ \ h(P, f(E))=h(f(P),f(E))=h(P,E)
\end{equation}
because $P$ is invariant under $f$ for any $f\in H$. Then, for example,
$$\sum_{i=1}^{m}d(P_1,f(F_i))\leq c.$$
Taking the sum as $f$ runs through the elements of $H$ we have that $c'\leq 2pc$. On the other hand if $X$ is an outer point of $f(E)$ for any $f\in H$ then
$$\sum_{i=1}^{m}d(X,f(F_i))>c \ \ (f\in H),\ \ \textrm{i.e.}\ \ \sum_{f \in H}\sum_{i=1}^{m}d(X,f(F_i))>2pc\geq c'.$$
By contraposition, if $X$ belongs to the polyellipse $E_{sym}$ then it belongs to the convex hull of $f(E)$ for {\emph{some}} $f\in H$. Finally
$$E_{sym}\subset \bigcup_{f \in H}\ \textrm{conv}\  f(E),$$
where  $\textrm{conv}\  f(E)$ denotes the convex hull of the polyellipse $f(E)$ $(f\in H)$. Since each polyellipse is around $P$ it follows that 
$$h(P, E_{sym})\leq h \left(P, \bigcup_{f \in H}\ \textrm{conv}\  f(E)\right)\leq h\left(P, \bigcup_{f \in H}\ f(E)\right)=h(P,E)$$
because the possible distances between the points of the sets $P$ and $f(E)$ are independent of the choice of $f\in H$ (see the minimax characterization (\ref{minimax}) of the Hausdorff distance). Figure 2 (left hand side) illustrates an ellipse $E$ with focuses $F_1$, $F_2$ and its symmetric pairs with respect to the isometry group of the triangle ($p=3$). The focal set of $E_{\textrm{sym}}$ (right hand side) contains three different points $F_1, F_2, F_3$ (solidcircles) with multiplicity $k_1=k_2=k_3=2$. The constant is choosen such that $E_{\textrm{sym}}$ passes through the vertices of the polygon.

\begin{figure}
\centering
\includegraphics[scale=0.3]{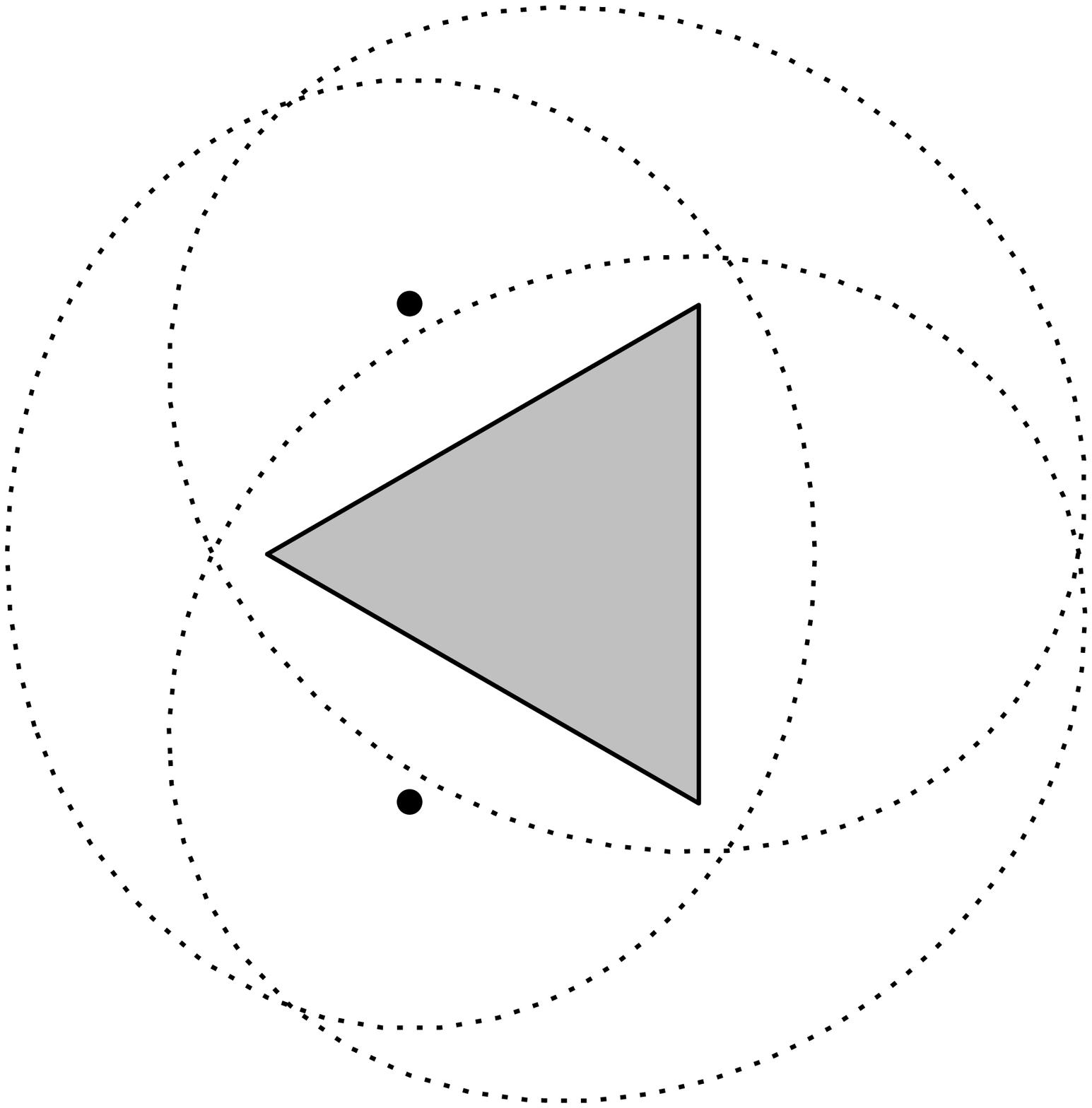}
\includegraphics[scale=0.3]{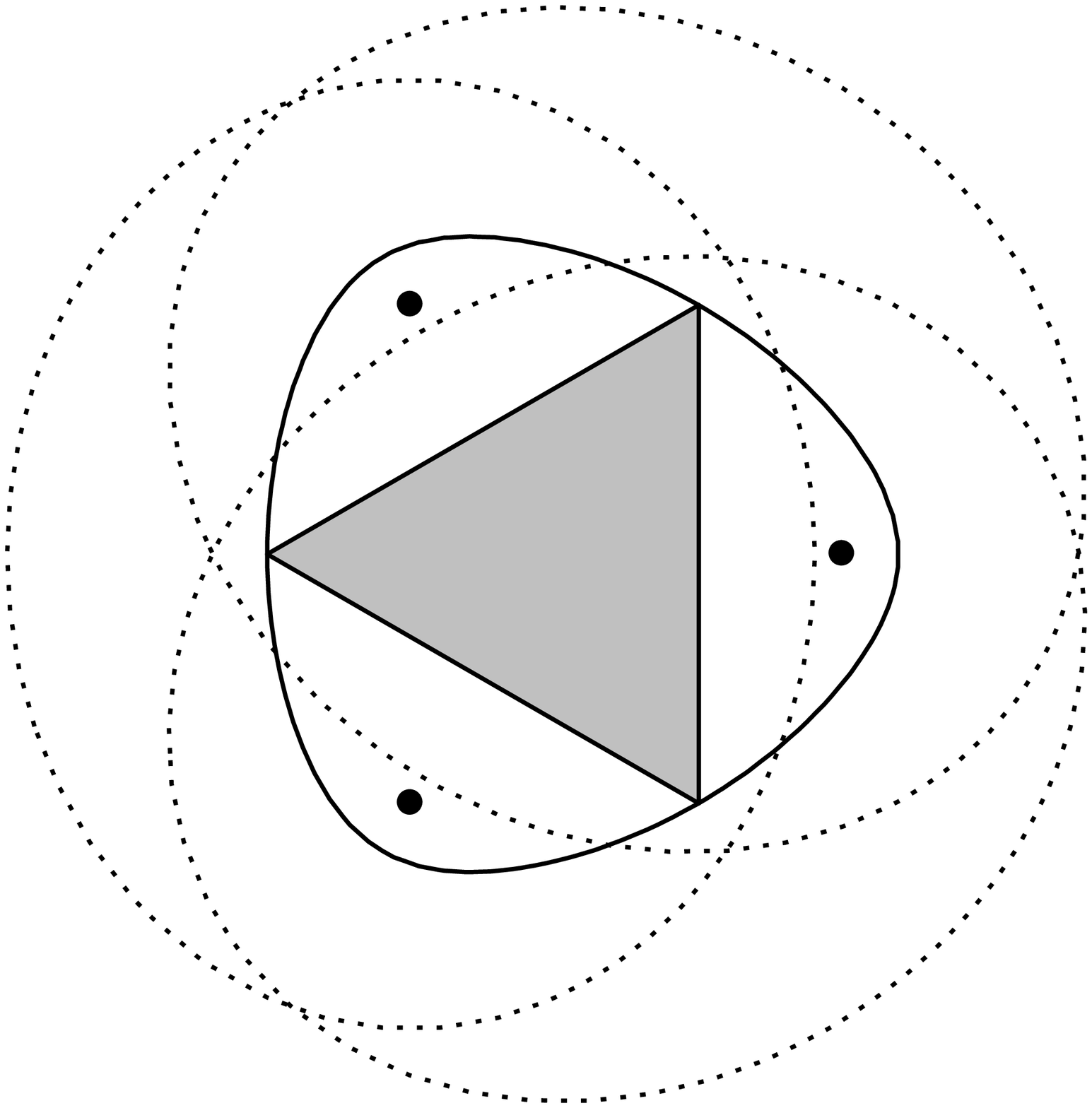}
\caption{The symmetrization process}
\end{figure}

\subsubsection{Summary} From now on we suppose that $E_n$ is a sequence of circumscribed polyellipses tending to $P$ such that each element of the sequence passes through the vertices $P_1, \ldots, P_p$ of $P$ and the set of the focuses is of the form (\ref{symmetry}).

\subsection{The third step - smoothness.} By the symmetry properties it can be easily seen that
$$h(P, E)=d(M_i, Q_i)\ \ (i=1, \ldots p)$$
where $E$ stands for a general element of the sequence of polyellipses, $M_i$ and $Q_i$ are the corresponding midpoints of the edge\footnote{The indices are taken by modulo $p$.} $P_iP_{i+1}$ of $P$ and the arc $P_iP_{i+1}$ of the polyellipse, respectively - note that the line passing through $Q_i$ parallel to the edge $P_iP_{i+1}$ of the polygon supports the polyellipse at $Q_i$ because of the symmetry, see Figure 3. For the sake of simplicity consider the case of $i=1$ and 
$$Q:=Q_1, \ \ M:=M_1.$$
It is clear that the polyellipse $E$ heritages the symmetries of its focal set, i.e. $E$ is invariant under the reflection about the horizontal line. So are the supporting lines at $Q$. Therefore the line passing through $Q$ parallel to the edge $P_1P_{2}$ of the polygon supports the polyellipse at $Q$. 

Since the next step of the proof will be the estimation of the curvature of the polyellipse $E$ at the point $Q$ but a polyellipse can pass through some of its focuses (see Figure 1), we need an extra process to avoid singularities. Consider the equation   
$$\sum_{f \in H}\sum_{i=1}^{m}d(X,f(F_i))=c$$
of the polyellipse $E$, where 
$$c:=\sum_{f \in H}\sum_{i=1}^{m}d(P_1,f(F_i))=\sum_{f \in H}\sum_{i=1}^{m}d(P_2,f(F_i))= \ldots =\sum_{f \in H}\sum_{i=1}^{m}d(P_p,f(F_i))$$
and suppose that $F_1=Q$ such that
$$\sum_{f \in H}\sum_{i=1}^{m}d(Q,f(F_i))=c.$$
Then 
$$\sum_{f \in H}\sum_{F_i=Q} d(Q,f(Q))+\sum_{f \in H}\sum_{F_i\neq Q}d(Q,f(F_i))=c,$$
\begin{equation}
\label{Q}
k_1\sum_{f \in H}d(Q,f(Q))+\sum_{f \in H}\sum_{F_i\neq Q}d(Q,f(F_i))=c
\end{equation}
and, because of $P_1\in E$, 
\begin{equation}
\label{q}
k_1\sum_{f \in H} d(P_1,f(Q))+\sum_{f \in H}\sum_{F_i\neq Q}d(P_1,f(F_i))=c.
\end{equation}
\begin{figure}
\centering
\includegraphics[scale=0.28]{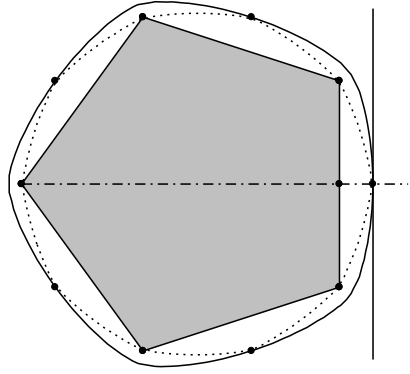}
\caption{Smoothness.}
\end{figure}
Consider the function
$$F(X):=\sum_{f \in H} d(X,f(M));$$
using the triangle inequality it can be easily seen that $F$ is a strictly convex function such that
$$F(P_1)=\ldots=F(P_n).$$
Therefore $F(M) < F(P_1)$ because $M$ is the midpoint of the segment $P_1P_2$. If $Q$ is close enough to $M$ then
$$\sum_{f \in H} d(Q,f(Q))\approx F(M) < F(P_1) \approx \sum_{f \in H} d(P_1,f(Q))$$
by a continuity argument. Therefore, equations  (\ref{Q}) and (\ref{q}) show that  
$$c':=\sum_{f \in H}\sum_{F_i\neq Q}d(Q,f(F_i)) > \sum_{f \in H}\sum_{F_i\neq Q}d(P_1,f(F_i)).$$
Now the polyellipse $E'$ defined by
$$\sum_{f \in H}\sum_{F_i\neq Q}d(X,f(F_i))=c'$$
passes through $Q$ but $P_1$ and, by the symmetry, $P_2$ are in the interior of $E'$ (see Figure 3). Since $Q$ (together with all symmetric pairs) has been deleted from the focal points we can take the curvature $\kappa$ of $E'$ at $Q$:
$$\kappa (Q) \leq \ \frac{1}{\textrm{the radius of the circle passing through $P_1, P_2$ and $Q$}}.$$
In case of a sequence of polyellipses tending to $P$, the point $Q$ can be arbitrarily close to $M$ and $\kappa(Q)$ can be arbitrarily close to zero.

\subsection{The estimation of the curvature}

From now on we suppose that the set of the focuses of the polyellipse $E$ is of the form (\ref{symmetry}) and $Q$
is the corresponding point on the arc of the polyellipse to the midpoint $M$ of the edge $P_1P_2$. To simplify the curvature formula \cite{G} at $Q$ as far as possible suppose that 
$$P_1(\cos(\pi/p),-\sin(\pi/p)),\ P_2(\cos(\pi/p),\sin(\pi/p))\ \ \textrm{and}\ \ M(\cos(\pi/p),0),$$
i.e. the tangent line to the polyellipse at $Q$ is parallel to the $y$ - axis. Therefore $D_2F(Q)=0,$
where
\begin{equation}
\label{curvformula}
F(X)=\sum_{f \in H}\sum_{i=1}^{m}d(X,f(F_i))\ \ \textrm{and}\ \ \kappa(Q)=\frac{D_2D_2F(Q)}{D_1F(Q)}.
\end{equation}
Since $D_1F(Q)$ is the first coordinate of the gradient vector of $F$ at $Q$ we have from (\ref{firstder}) that
$$0 < D_1F(Q)=\| \sum_{f \in H}\sum_{i=1}^{m}\frac{1}{d(Q,f(F_i))}\stackrel{\longrightarrow}{Qf(F_i)}\|\leq 2pk_0+ \| \sum_{f \in H}\sum_{F_i \neq O}\frac{1}{d(Q,f(F_i))}\stackrel{\longrightarrow}{Qf(F_i)} \|,$$
where $k_0$ is the multiplicity of the symmetry center $O$ if it belongs to the set of focuses $F_1$, $\ldots$, $F_m$ and $k_0=0$ otherwise. Applying Theorem 1 to the function 
$$F_0(X):=\sum_{f \in H}\sum_{F_i \neq O}d(X,f(F_i))$$
it follows that 
$$\sum_{f \in H}\sum_{F_i \neq O}\frac{1}{d(O,f(F_i))}\stackrel{\longrightarrow}{O f(F_i)}={\bf 0}.$$
Therefore
$$0 < D_1F(Q)=\| \sum_{f \in H}\sum_{i=1}^{n}\frac{1}{d(Q,f(F_i))}\stackrel{\longrightarrow}{Qf(F_i)}\|\leq 2pk_0+ \| \sum_{f \in H}\sum_{F_i \neq O}\frac{1}{d(Q,f(F_i))}\stackrel{\longrightarrow}{Qf(F_i)} \|= $$
$$2pk_0+\| \sum_{f \in H}\sum_{F_i \neq O}\frac{1}{d(Q,f(F_i))}\stackrel{\longrightarrow}{Qf(F_i)} - \sum_{f \in H}\sum_{F_i \neq O}\frac{1}{d(O,f(F_i))}\stackrel{\longrightarrow}{Of(F_i)}\|\leq$$
$$2pk_0+\sum_{f \in H}\sum_{F_i \neq O} \| \frac{1}{d(Q,f(F_i))}\stackrel{\longrightarrow}{Qf(F_i)} - \frac{1}{d(O,f(F_i))}\stackrel{\longrightarrow}{Of(F_i)}\|.$$

\begin{lemma}
If $Q$ is close enough to $M$ then 
$$\| \frac{1}{d(Q,f(F_i))}\stackrel{\longrightarrow}{Qf(F_i)} - \frac{1}{d(O,f(F_i))}\stackrel{\longrightarrow}{Of(F_i)}\|\leq \frac{4}{1+d(O,F_i)}$$
for any $F_i\neq O$, $i=1, \ldots, m$.
\end{lemma}

\begin{Pf} If $F_i$ is in the interior of the unit circle centered at $O$ then the same is true for any point of the form $f(F_i)$ because $O$ is the fixpoint of the elements in the symmetry group $H$. Therefore
$$\frac{4}{1+d(O,F_i)}\geq \frac{4}{1+1}=2$$
which is the maximum of the length of the difference of unit vectors. Otherwise, if $F_i$ (together with its all symmetric pairs) is an outer point of the unit circle centered at $O$ and $Q$ is in its interior, i.e. $Q$ is close enough to $M$, then $\alpha_i:=\angle(Of(F_i)Q) < 90^{\circ}$ for any $f\in H$. By the cosine rule 
$$\| \frac{1}{d(Q,f(F_i))}\stackrel{\longrightarrow}{Qf(F_i)} - \frac{1}{d(O,f(F_i))}\stackrel{\longrightarrow}{Of(F_i)}\|^2=2(1-\cos \alpha_i)=4\sin^2 \frac{\alpha_i}{2},$$
i.e.
$$\| \frac{1}{d(Q,f(F_i))}\stackrel{\longrightarrow}{Qf(F_i)} - \frac{1}{d(O,f(F_i))}\stackrel{\longrightarrow}{Of(F_i)}\|=2\sin \frac{\alpha_i}{2}\leq 2\sin \alpha_i$$
because of $0\leq \alpha_i < 90^{\circ}$. By the sine rule
$$d(O,F_i)\sin \alpha_i \leq d(O,Q) < 1\ \ \Rightarrow\ \ \sin \alpha_i < \frac{2}{1+d(O,F_i)}$$
and the proof is complete. 
\end{Pf}

\begin{corollary}
$$0 < D_1F(Q_1)\leq 2pk_0+8p\sum_{F_i\neq O}\frac{1}{1+d(O,F_i)}\leq 8p\sum_{i=1}^m \frac{1}{1+d(O,F_i)}.$$
\end{corollary}

\begin{lemma} If $Q$ is close enough to $M$ then 
$$D_2D_2(Q) \geq 2\cos^2 \frac{\pi}{p} \sum_{i=1}^{m}\frac{1}{1+d(O,F_i)},$$
\end{lemma}

\begin{Pf}
Using formula (\ref{secondder}) we should estimate expressions of type 
\begin{equation}
\label{masod}
S(Q,F_i):= \sum_{f \in H}\frac{1}{d(Q,f(F_i))}\cos^2 \beta_i
\end{equation}
for any $i=1, \ldots, m$, where $\beta_i=\angle (F_iQO)$. By the triangle inequality 
$$d(Q,f(F_i)) \leq d(Q,O)+d(O,f(F_i)) \leq 1+d(O,F_i)$$
provided that $Q$ is in the interior of the unit circle centered at $O$, i.e. $Q$ is close enough to $M$. Therefore
\begin{equation}
\label{estimation1}
\frac{1}{1+d(O,F_i)}\sum_{f \in H}\cos^2 \beta_i \leq S(Q,F_i).
\end{equation}
First of all consider the case of $F_i=O$: by (\ref{estimation1})
\begin{equation}
\label{spec}
S(Q,O)\geq \sum_{f \in H} 1=2p \geq 2\cos^2 \frac{\pi}{p}\frac{1}{1+d(O,O)}.
\end{equation}
Suppose that $F_i\neq O$. Since the focal set contains all symmetric pairs of $F_i$ by the action of the elements in $H$ there must be at least two focal points of the form $f(F_i)$ in the sector of the plane with polar angle between $\pi-\pi/p$ and $\pi+\pi/p$ (one of them as a rotated and another one as a reflected pair of $F_i$), see Figure 4.  
\begin{figure}
\centering
\includegraphics[scale=0.3]{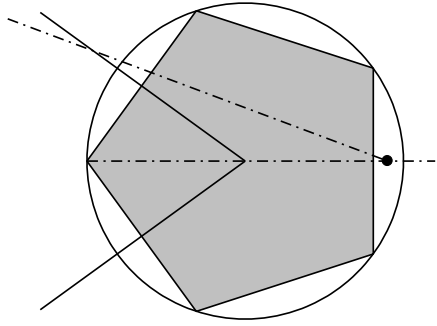}
\caption{The proof of Lemma 2.}
\end{figure}
Therefore
\begin{equation}
\label{alt}
S(Q,F_i) \geq 
\frac{1}{1+d(O,F_i)}\sum_{f \in H}\cos^2 \beta_i \geq 2\cos ^2 \frac{\pi}{p}\frac{1}{1+d(O,F_i)}
\end{equation}
as was to be proved.
\end{Pf}

Using the curvature formula (\ref{curvformula}) Corollary 1 and Lemma 2 gives that
\begin{equation}
\label{becsles}
\kappa(Q) \geq \frac{\cos^2 \frac{\pi}{p}}{4p}
\end{equation}
independently of the number of focuses. Therefore $\kappa(Q)$ can not be arbitrarily close to zero which is a contradiction. 

\section{Concluding remarks} To present a reach class of curves in the plane as Hausdorff limits of polyellipses we should use integration to compute the average distance of a point from a given (focal) set. Let $\Gamma$ be a plane curve with finite arclength $L(\Gamma)$ and consider the function
\begin{equation}
\label{gencon}
f(X):=\frac{1}{L(\Gamma)}\int_\Gamma d(X,Y)\, dY.
\end{equation}
The curve given by the equation $f(X)=c$ is called a generalized conic with focal set $\Gamma$. In what follows we show that such a curve is the Hausdorff limit of polyellipses. First of all note that $f$ is convex because of the convexity of the integrand. This means that it is a continuous function, i.e. the level set $f^{-1}(c)$ is closed. On the other hand it is bounded because $\Gamma$ is bounded. Therefore a generalized conic is a compact subset in the plane.  Let $\varepsilon >0$ be a given positive number and consider the partition of $\Gamma$ into $2^m$ equal parts (depending on $X$) such that 
\begin{equation}
\label{upperlower}
0<S(X,m)-s(X,m)< L(\Gamma) \frac{\varepsilon}{2},
\end{equation}
where $X\in f^{-1}(c)$, $S(X,m)$ and $s(X,m)$ are the upper and lower Riemann sum of
$$\int_\Gamma d(X,Y)\, dY$$
belonging to the equidistant partition. By a continuity argument\footnote{Recall that both the sup-function and the inf-function are Lipschitzian:
$$\inf_{\xi} d(X, \xi)\leq \inf_{\xi} d(X,Y)+d(Y, \xi)=d(X,Y)+\inf_{\xi} d(Y, \xi),$$
i.e. $|\inf_{\xi} d(X, \xi)-\inf_{\xi} d(Y, \xi)|\leq d(X,Y).$ In a similar way, 
$$\sup_{\xi} d(X, \xi)\leq \sup_{\xi} d(X,Y)+d(Y, \xi)=d(X,Y)+\sup_{\xi} d(Y, \xi).$$} formula (\ref{upperlower}) holds on an open disk centered $X$ with radius $r_X$. Since $f^{-1}(c)$ is compact we can find a finite open covering with open disks centered at some points $X_1$, $\ldots$, $X_k\in f^{-1}(c)$. Let
$$m:=\max \{m_1, \ldots, m_k\},$$ 
where $m_i$ denotes the power of the partition of $\Gamma$ into $2^{m_i}$ equal parts. Obviously we have a refinement of all partitions under the choice of the maximal value of $m_i$'s. Therefore
\begin{equation}
\label{upperlowerglobal}
0<S(X,M)-s(X,M)< L(\Gamma) \frac{\varepsilon}{2}
\end{equation}
for any $X\in f^{-1}(c)$, where $M:=2^{m}$. Let $\tau_1$, $\ldots$, $\tau_M\in \Gamma$ be some middle points of the partition and define the function
\begin{equation}
\label{approx}
F_M(Z):=\frac{d(Z,\tau_1)+\ldots+d(Z,\tau_M)}{M}.
\end{equation}
Since we have an equidistant partition it follows that
\begin{equation}
\frac{1}{L(\Gamma)}s(X,M) \leq F_M(X) \leq \frac{1}{L(\Gamma)}S(X,M) \ \ \Rightarrow \ \ |F_M(X)-c|< \varepsilon
\end{equation}
for any $X\in f^{-1}(c)$. This means that $f^{-1}(c)$ is between the polyellipses defined by the equations $F_M(X)=c-\varepsilon$ and $F_M(X)=c+\varepsilon$, i.e. $\{ X | F_M(X)=c-\varepsilon\} \ \subset \ f^{-1}(c)\ \subset \{ X | F_M(X)=c+\varepsilon\}.$ Taking the limit $\varepsilon \to 0$ the sequences of the polyellipses on both the left and the right hand sides of the previous formula tends to $f^{-1}(c)$. 

\subsection{A Maple presentation II}

Since it is hard to find the equidistant partition of a curve in general we use the equidistant partition of the parameter interval to present an explicite example: let $\Gamma: \ \gamma\colon [0,2\pi]\to \gamma(t)=(t,\sin(t))$ be a period of the sine wave. Figures 5 and 6 show the generalized conic $f^{-1}(4)$ in "dot" linestyle (see equation (\ref{gencon})). The approximating polyellipses belong to the equidistant partition of the interval into 
$M=2^3$, $2^4$, $2^5$ and $2^6$ equal parts, respectively. 

\begin{figure}
\centering
\includegraphics[scale=0.25]{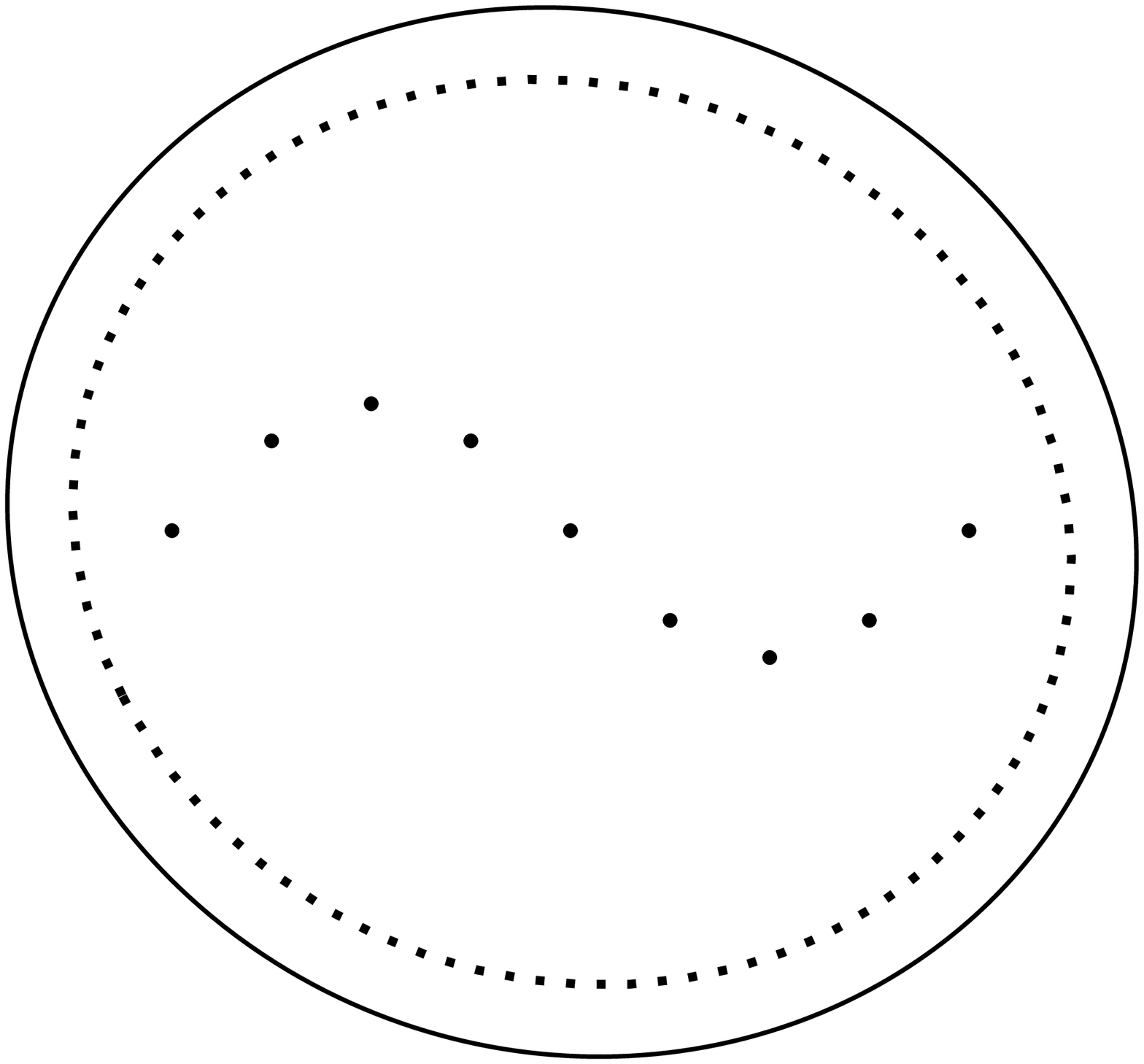}
\includegraphics[scale=0.25]{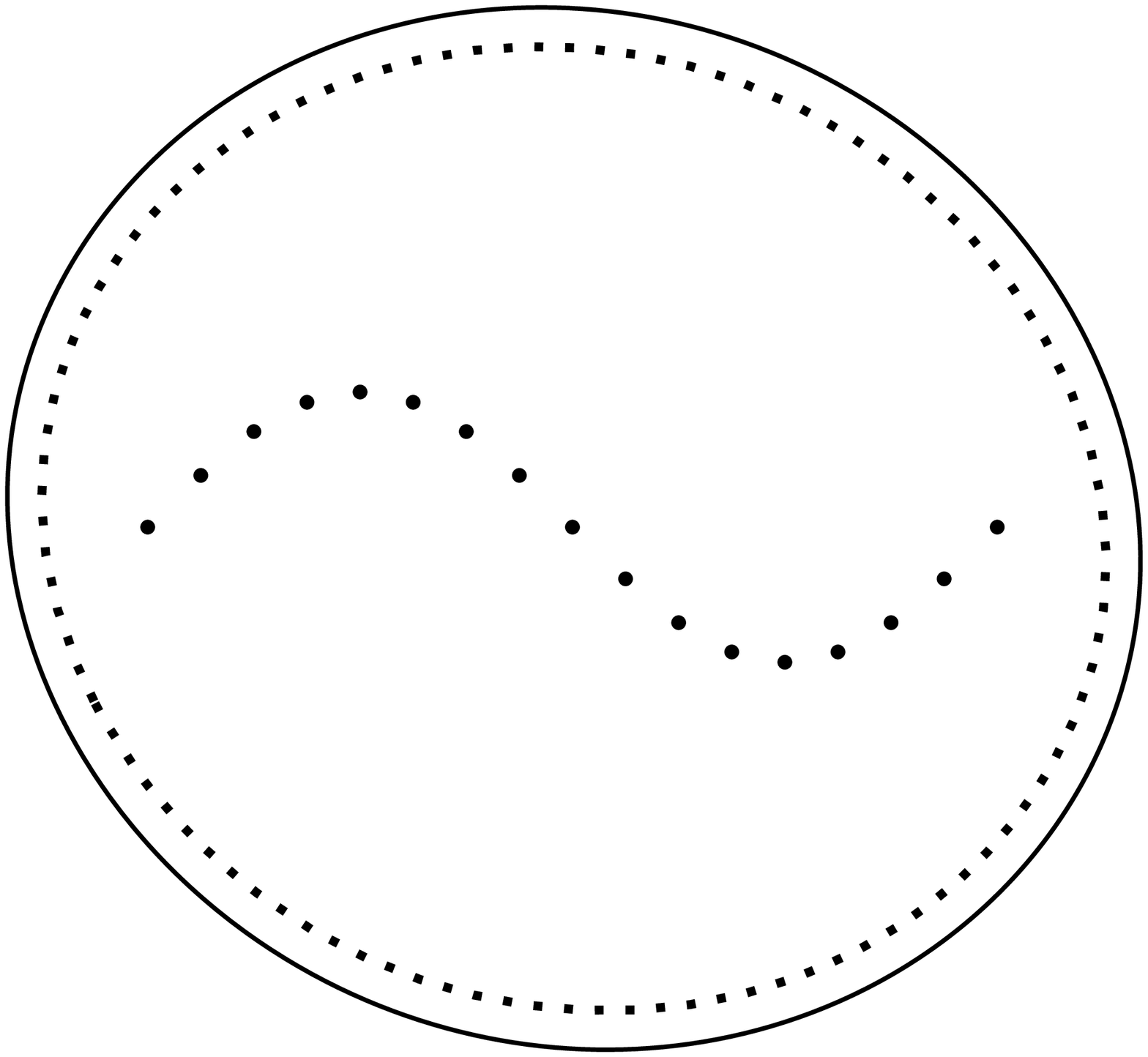}
\caption{A Maple presentation II: the case of $M=2^3$ and $2^4$.}
\end{figure}

\begin{verbatim}
with(plottools):
with(plots):
#The coordinates of the focal points#
x:=[seq((2*Pi*(1/8))*k,k=0..8)]:
y:=[seq(sin((2*Pi*(1/8))*k),k=0..8)]:
m:=nops(x):
F:=proc(u,v)
local i,partial:
partial:=0.:
for i to m do 
partial:=partial+\frac{1}{m}sqrt((u-x[i])^2+(v-y[i])^2) 
end do 
end proc:
#The list of the level rates#
c:=[4]:
graf1:=pointplot(zip((s,t)->[s,t],x,y),symbol=solidcircle,
symbolsize=10, color=black):
graf2:=seq(implicitplot(F(u,v)=i,u=-2..8,v=-5..5,
numpoints=10000,color=black),i=c):
#Integration along the sine wave#
F:=(u,v)->(1/4)*Int(sqrt((u-t)^2+(v-sin(t))^2)*
sqrt(1+cos(t)^2),t=0..2*Pi)/(sqrt(2)*EllipticE((1/2)*sqrt(2))):
graf3:=implicitplot(F(u,v)=4,u=-2..8,v=-5..5,linestyle=dot,
thickness=3,color=black):
display(graf1,graf2,graf3,scaling=constrained,axes=none);
\end{verbatim}

\begin{figure}
\centering
\includegraphics[scale=0.25]{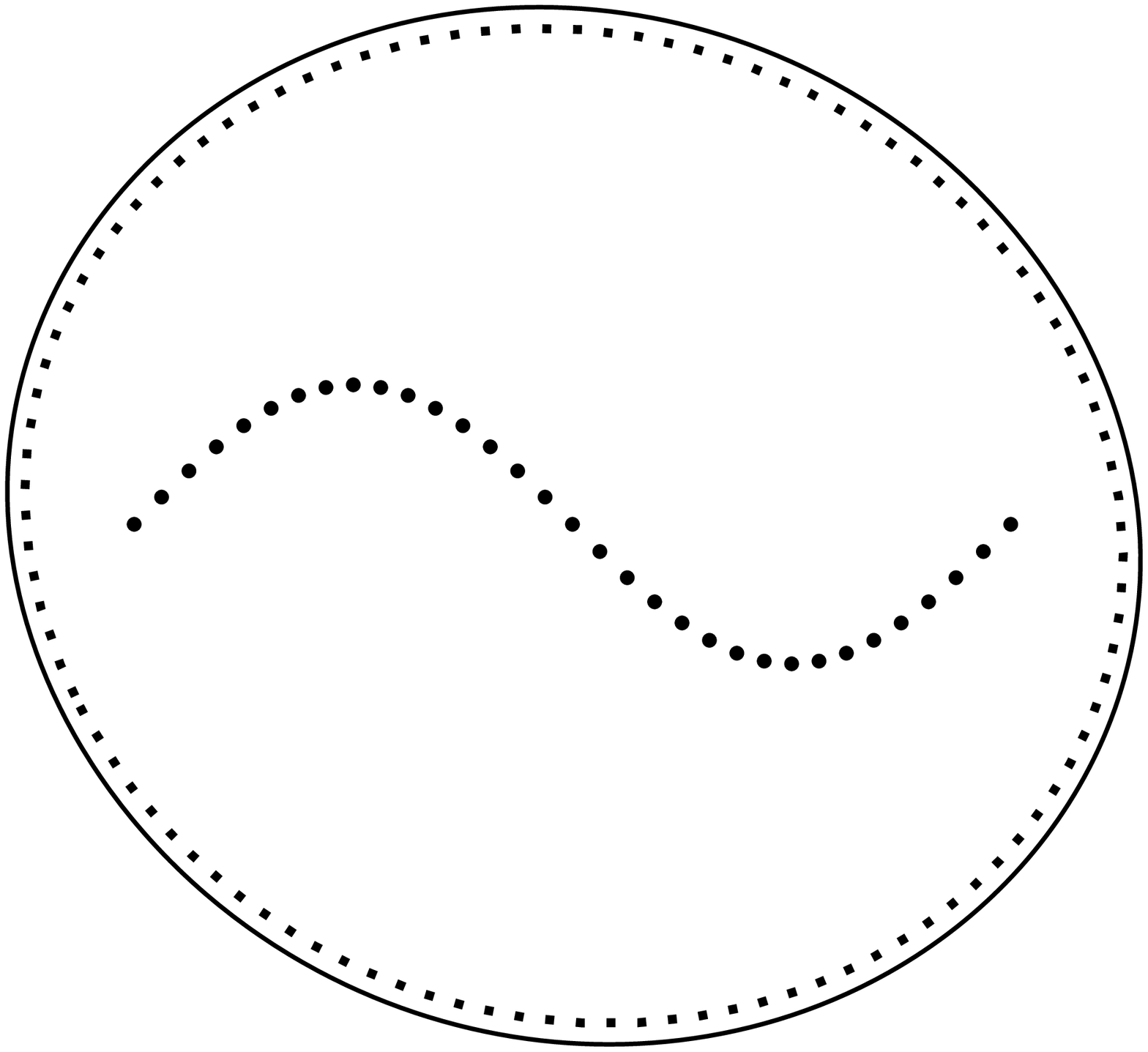}
\includegraphics[scale=0.25]{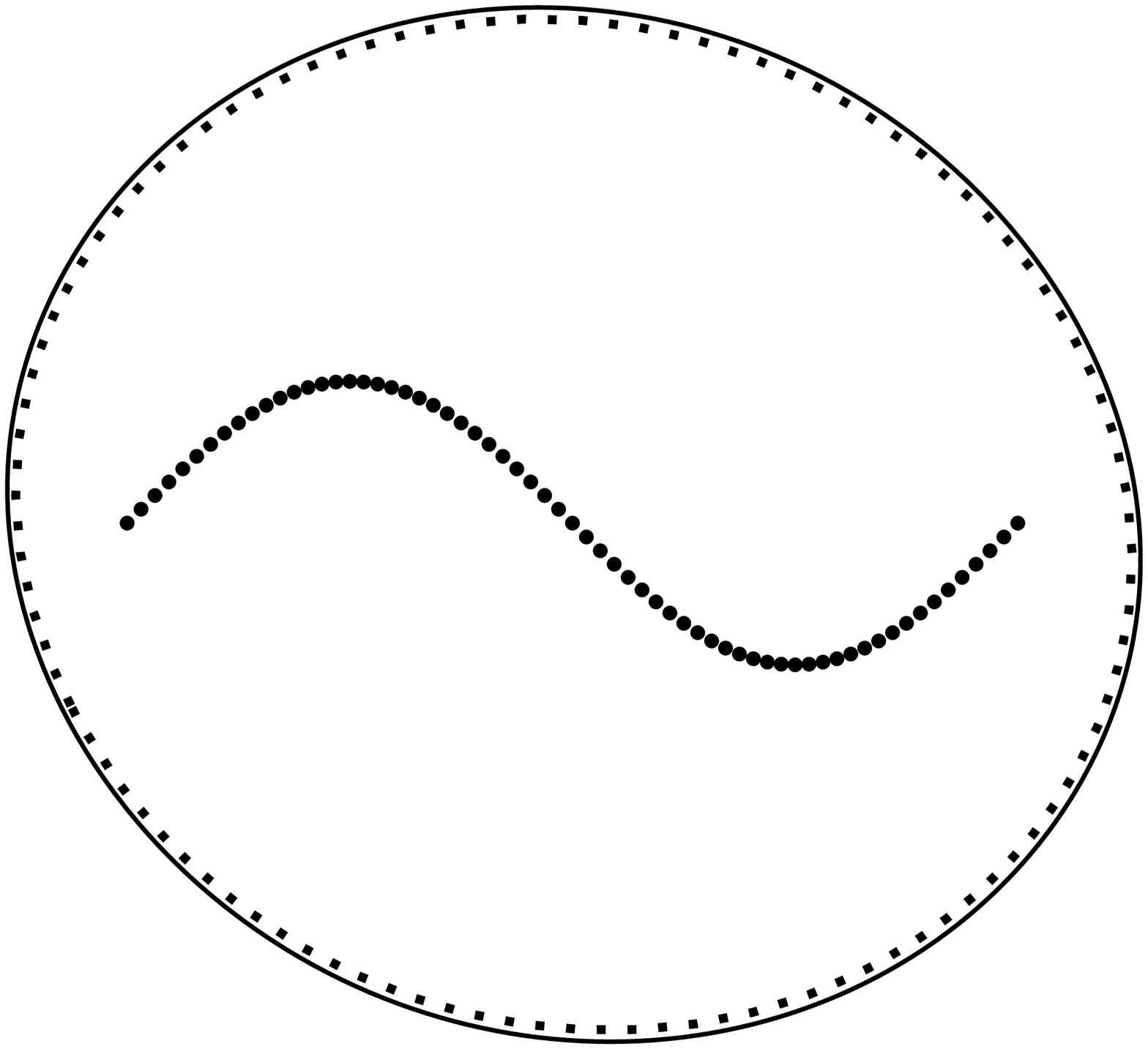}
\caption{A Maple presentation II: the case of $M=2^5$ and $2^6$.}
\end{figure}

\subsection{A note about the error of the approximation of a regular polygon by polyellipses}
If $\Gamma$ is the unit circle in the plane then $f^{-1}(c)$ is invariant under any element of the isometry group because the level sets heritage the symmetries of the integration domain (focal set). Therefore $f^{-1}(c)$ is a circle centered at the same point as $\Gamma$. Especially, $\Gamma=f^{-1}(c)$, where $c=4/\pi$. This means that the error of the approximation of a regular polygon by polyellipses tends to zero as the number of the vertices tends to the infinity: if $P$ is a regular $p$ - gon inscibed in the unit circle $\Gamma$, then 
$$h(P,\Gamma)=1-\cos\frac{\pi}{p}\ \ \textrm{and}\ \ h(P,E)\leq h(P, \Gamma)+h(\Gamma, E)\leq 1-\cos\frac{\pi}{p}+\varepsilon$$  
provided that $E$ is an approximating polyellipse of $\Gamma$ with error less then $\varepsilon$.

\end{document}